\documentclass[a4paper, 10pt]{article}
\usepackage{mathrsfs}
\usepackage{amsfonts}
\usepackage{latexsym}
\usepackage{graphicx}
\usepackage{amsmath}
\newcommand{\new}{\newcommand*}

\new{\rnew}{\renewcommand*}
\new{\newe}{\newenvironment*}
\new{\newt}{\newtheorem}
\new{\stl}{\setlength}

\new{\bea}{\begin{eqnarray}}
\new{\eea}{\end{eqnarray}}

\new{\be}{\begin{equation}}
\new{\ee}{\end{equation}}

\new{\bean}{\begin{eqnarray*}}
\new{\eean}{\end{eqnarray*}}

\new{\no}{\nonumber}

\new{\bt}{\begin{theorem}}
\new{\et}{\end{theorem}}
\new{\bl}{\begin{lemma}}
\new{\el}{\end{lemma}}

\new{\bc}{\begin{corollary}}
\new{\ec}{\end{corollary}}
\new{\bp}{\begin{proof}\quad}
\new{\ep}{\end{proof}}
\new{\ba}{\begin{array}}
\new{\ea}{\end{array}}

\newt{prop}{Proposition}
\newt{lem}{Lemma}
\newt{conjecture}{Conjecture}[section]
\newt{cor}{Corollary}
\newt{thm}{Theorem}
\newt{assumption}{Assumption}
\newt{ex}{Example}
\newt{remark}{Remark}
\newt{definition}{Definition}

\rnew{\theequation}{\thesection.\arabic{equation}}
\new{\sect}[1]{ \section{#1}
\setcounter{equation}{0} \setcounter{figure}{0} }

\newe{acknowledgment}{{\flushleft \bf Acknowledgments. }}{}
\newe{proof}{{\bf Proof.}}{\qquad\endproof}
\newe{keywords}{\begin{quote}\quad{\bf Keywords.}}{\end{quote}}
\newe{AMS}{\begin{quote}\quad {\bf AMS subject classification.}}{\end{quote}}

\def \endproof {\qquad \vrule height 5pt width 5pt depth 2pt}

\title{Continuity of the inverse of   Riemann map\date{}}

\author{{Zhijian Qiu  }
\\{\small   School of  mathematics}
\\{\small  Southwestern University of Finance and Economics}
\\{\small   Chengdu,  China, E-mail:\ qiu@swufe.edu.cn }
}

\begin{document}
\maketitle
\begin{abstract}
For a simply connected domain  $G$,  let $\partial_{a}G$ be the set
of $\lambda$ in $\partial G$ for which   there is a Jordan arc
$J\subset \overline G$ such that $  J\cap\partial G=\{\lambda\}$.
Let $\partial_{n} G=\partial G-\partial_{a}G$. A point
$\lambda\in\partial G$ is called semi-unreachable  if there is a
simply connected subdomain  $ U$ such that $\lambda\in
\partial_{n}U$ and $\lambda\ \overline \in\ \overline{\partial U\cap
G}$. Let  $\partial_{sn}G$ denote the set of semi-unreachable
points. We show that  the reverse of Riemann map
  from the unit disk $D$ onto $G$  extends
continuously to $\overline D$ if and only if
$\partial_{sn}G=\emptyset  $.

As a consequence, we obtain a very short and elementary proof for
the Osgood conjecture: if $G$ is a Jordan domain, then 
the Riemann map of $G$ extends to be a homeomorphism from $\overline
G$ to $\overline D$.\footnote{ This is a very famous result and   is
known
as the Carathe\"{o}dory theorem.}
\end{abstract}
\begin{keywords} Riemann map,  complex analysis, accessible point.
\end{keywords}
\begin{AMS} 30H05, 30E10.
\end{AMS}

 \section*{Introduction}

 The Riemann mapping theorem is one of most important results in classical mathematics. It states that a simply connected
domain $G$ whose boundary does not reduce  to a single point can be
mapped onto the unit disk $D$ by a univalent analytic function
 $\varphi$. Poincar\'{e} showed that $\varphi$ is essentially unique  and it is called Riemann  map.
In this article we study the problem of when $\psi$, the inverse of
$\varphi$, extends continuously to  $\partial D$? The boundary
behavior of both $\varphi$ and $\psi$ were investigated intensively
by many authors throughout history.  Schwarz and Painlev\'{e} and
others  proved that $\psi$ is a homeomorphism from $\overline D$
onto $\overline G$ if $G$ is a domain with piecewisely smooth
boundary.  Osgood   conjectured the result  would be true for  a
Jordan domain  in 1900,  and
  Carathe\"{o}dory proved Osgood's conjecture in 1913 \cite{cara1913}.

 In the meanwhile, Carathe\"{o}dory developed the theory of prime ends in \cite{cara1913b},
   a general and vast theory   for the boundary behavior of the Riemann map from a general
simply connected domain to    D. It can be viewed as a kind of
extension of the result for Jordan domains. This theory leads to the
result:
  $\psi$ extends continuously to $\partial D$ if and only if $\partial G$ is locally connected (\cite{pomm},
p.20).  However, it seems that the majority of people in the math
society do not know this result.

 The purpose of this article is to seek a
direct, elementary  and classical solution for the problem.
 The motivation is   the strong desire for a theorem readable for  someone who just had a
graduate course in complex analysis. We begin  by introducing the
key concept: semi-unreachable   points, which is purely topological
and understandable to students. It turns out that with such a
concept we are able to provide  a very simple and elementary
solution: $\psi$ extends continuously to $\partial D$ if and only if
$\partial G$ has no semi-unreachable points. One advantage of our
theorem is that it not only offers an easy and direct solution
understandable to ordinary  mathematicians but also leads to a short
and concise proof for the Osgood conjecture. Our proof uses only
some basic results established already a century ago and even does
not rely on the Lebesgue theory. Over the past few decades, the
author was unable to find a proof for the Carathe\"{o}dory theorem
in literatures understandable to the general math public. In fact,
L. Ahlfors wrote in \cite{a}[p.232]: "Unfortunately, considerations
of space do not permit us to include a proof of this important
theorem (the proof would require a considerable amount of
preparation)". The theorem he referred to was the Carathe\"{o}dory
theorem. Actually, it seems that the Carathe\"{o}dory theorem
appears only in literatures for the specialized fields. After a long
time of efforts, the author now is able to provide this work which
gives a solution understandable to the math public and  can be
easily presented in standard text books like \cite{a} and \cite{b}
(since they cover all the materials needed for this article). A
desirable solution should be teachable along with the Riemann
mapping theorem in a standard course. In this respect, our theorem
is the first desirable solution most close to the essence
of the original problem.

Lastly, our method is localizable to determine on which part a
univalent analytic function   continuously extends to the boundary
for a general domain.

\section{The Result}

A Jordan curve $\gamma$  in the complex plane is   the image of the
interval $[a,b]$  under a continuous function $f$ such that $f$ is
injective on $(a,b)$ and $f(a)=f(b)$. $\gamma $ is called a Jordan
arc if it the image of  an injective continuous function on $[a,b]$.
For a simply connected domain $G$, a crosscut  of $G$ is a Jordan
arc $J$ for which  the interior $J^{\circ}$ is contained in $ G$ and
$J\cap\partial G$ contains two the endpoints.
 For a domain  $G$, $\lambda\in
\partial G$ is said to be an accessible point  of $G$ if
there is a Jordan arc $J$ that is contained in $\overline G$ and
$J\cap\partial G =\{\lambda\}$. Let $\partial _{a}G$ be the set of
accessible points  and let $\partial_{n} G=\partial
G-\partial_{a}G$.

We use  $\varphi$ to denote the Riemann map and $\psi$ to denote the
reverse of   $\varphi$ throughout the paper.
The appendix contains the preliminaries for  this work.

\begin{definition} A point $\lambda\in\partial G$ is called semi-unreachable  if
there is a simply connected subdomain  $ U$ such that $\lambda\in
\partial_{n}U$ and $\lambda\ \overline \in\ \overline{\partial U\cap G}$.
\end{definition}

\begin{prop} \label{p}
$\lambda\in \partial_{sn} G$ if and only if there is a crosscut
$\gamma$ of $G$ and domains $U_{1}$ and $U_{2}$ such that
$G-\gamma=U_{1}\cup U_{2}$ and $\lambda\in(\partial_{n}
U_{1}\bigcup\partial_{n} U_{2})-\gamma$.
\end{prop}
\begin{proof}
Sufficiency  is straightforward  from the definition, so we only
prove the necessity. Let $\lambda\in\partial _{sn}G$, then there is
$U\subset G$ such that $\lambda\in\partial_{n}U$ and $\lambda\
\overline \in\ \overline{\partial U\cap G}$. Let $W=\psi^{-1}(U)$,
then there is a component of $l$ of $\partial W\cap D$ such that
$D-l$ has at least two components, say $W_{1}$ $\&$ $W_{2}$, and for
which $W\subset W_{1}$, $
\partial W_{2}\cap\partial D \neq\emptyset$ and $(\partial
W_{2}\cap\partial D) ^{\circ}\cap\partial W=\emptyset$. Let $\gamma$
be a crosscut of $G$ that is  contained in $\psi(W_{2})$ and let
$G-\gamma=U_{1}\cup U_{2}$, where $U_{1} $ is the component that
contains $U$ and
$U_{2}$ is the other one. 
We now show that $\lambda\in\partial_{n} U_{1}$. If not, there is a
Jordan arc $J\subset \overline U_{1}$ such that $\overline
J\cap\partial U_{1}=\{\lambda\}$. We may assume that  $diam(J)
<dist(\lambda,\partial U\cap G)$. Then the hypothesis implies that
$J^{\circ}\cap\partial U=\emptyset$. Note,  $dist(\lambda, \partial
U_{1}-\partial U)\geq dist(\lambda,\partial U\cap G)$, so
$J^{\circ}\cap(\partial U_{1}-\partial U)=\emptyset$, and it follows
that $J^{\circ}\subset U$. This means that
$\lambda\in\partial_{a}U$, contradicting the hypothesis. So
$\lambda\in\partial_{n}U_{1}$.
\end{proof}
\begin{thm} \label{main}     $\psi$ extends continuously to $\overline D$ if and only if
$\partial_{sn}G=\emptyset$.
\end{thm}
\begin{proof}
Necessity. 
Let  $a\in \partial G$ and let $\{z_{n}\}$ be a sequence in $  G$
such that $z_{n}\rightarrow a$.
 Then there is a subsequence $\{z_{n_{i}}\}$
    such that  $\{\varphi(z_{n_{i}})\}\rightarrow b$ for some
 $b\in\partial D$. There is a Jordan arc $J $ for which
 $J^{\circ}\subset D$ and
 $J$ contains a subsequence of $\{\varphi(z_{n_{i}})\}$.
 The hypothesis that
 $\psi$ is continous implies  that $\psi(J)$ is a Jordan arc  in $\overline  G$
 and $ \psi(J) \cap\partial G=\{a\}$. So
 $a\in\partial_{a}G$ and thus  $\partial G=\partial_{a} G$. Now if $\gamma$ is a crosscut
and $U$ $\&$ $V$ are domains such that $G-\gamma=U\cup V$, then
evidently no point in $\partial G-\gamma$ belongs to $
\partial_{n}U\cup
 \partial_{n} V$ (since $\partial U=\partial_{a}U$ and $\partial V=\partial_{a}V$).
 By Proposition~\ref{p},     $\partial_{sn}G =\emptyset$.

 Sufficiency.  If $\psi$ is
 not continuous, then there is $b\in\partial D$ and a sequence $\{z_{n}\}\subset D $
such that $z_{n}\rightarrow b$ but  $\{\psi(z_{n})\}  $ does not
converge. Then  there is a  Jordan arc $J$  such that
$J^{\circ}\subset D$, $J\cap\partial D=\{b\}$ and $J$ contains a
 subsequence of $\{z_{n}\}$ whose image under $\psi$ is divergent.  $\overline{ \psi(J)}\cap
\partial G$ is connected and is not a singleton. We now show that  of $\overline{ \psi(J)}\cap
\partial G \cap
\partial_{n}G\neq \emptyset$. In fact,   suppose the contrary. Then there are    $a$  \  $\&$\ $ b$
in $\overline{ \psi(J)}\cap \partial G\cap\partial_{a}G$, so
we have
 a crosscut of $G$ that connects $a$ and $b$. From our construction of $J$, we see
 that $\psi^ {-1}(a) =\psi^ {-1}(b)$, and the image of the crosscut under $\psi^{-1}$ would be a loop. By  Appendix E, it is
 a contradiction.
 So there is $y$ such that  $y\in \overline{ \psi(J)}\cap \partial
G\cap\partial_{n}G\neq
 \emptyset$. Let $U$ be a simply connected subdomain of $G$ such
 that the interior of  $\partial U \cap \partial G$ contains $\overline{ \psi(J)}\cap \partial
 G$, then $y\in \partial_{n}U$ and $y\ \overline \in\  \overline{\partial U\cap
 G}$ and  hence  $y\in \partial_{sn}G$, which contradicts the hypothesis.
 \end{proof} \\

We now have a direct and elementary  proof for the Osgood
conjecture.

\begin{cor} \label {cor}
If $G$ is a Jordan domain, then    $\psi$
  extends to be a homeomorphism. 
\end{cor}
\begin{proof}
If  $G$ is a Jordan domain, then for  any $\lambda\in
\partial G$ and any  $\delta >0$, there is a  Jordan domain $V\subset G$  with $diam
(V)<\delta$ so that $\partial U\cap\partial G$ is a subarc and
$\lambda  \in  (\partial U\cap\partial G)^{\circ}$. Now it is clear
that no subdomain $U$ satisfies the condition that $\lambda\in
\partial_{n}U$ and $\lambda\hspace{.02in} \overline \in\hspace{.05in} \overline{\partial U\cap G}$, and
so $\partial_{sn}G=\emptyset$. Therefore, $\psi$ extends
continuously to $\overline D$. For injectivity, suppose
$\psi(a)=\psi(b)$ for $a$ and $b$ in $\partial D$.  Let $J$ be a
crosscut that joins $a$ and $b$, then $\psi(J)$ is a Jordan curve
and $\psi(J)\cap\partial G=\{\psi(a)\}$. Let $U$ be the domain
enclosed by $\psi(J)$, then $\partial[ \psi^{-1}(U)]$ contains   a
subarc of
$\partial D$ 
and this implies  that $\psi$ is  constant  on $D$. It is a
contradiction and so  $\psi$ is a homeomorphism.
\end{proof}\\

The following results have been known for a century, however, the
 concise and elementary proofs below are due to the author's   constant efforts in years.\\

\noindent {\bf Appendix.}\\

 \noindent{\bf  A}). \
{\em A crosscut $J$ of   $ G$ separates $G$ into two disjoint simply
connected domains.}\\

{\bf Proof.} Let $G_{n}=\phi(\{z: |z|< 1-\frac{1}{n}\})$. Then for
all sufficiently large $n$, $G_{n}\cap J \neq \emptyset$.
Parameterize $J$ so that $ J(t): [0, 1] $ onto $J$. Let $s\in [0,1]$
such that $J(s) \in G_{n}$. Clearly $ J \cap
\partial G_{n} \neq \emptyset$. Let $a= \sup \{t: J(t) \in \partial G_{n}
\hspace{.1in}\mbox{and}\hspace{.1in} t\leq s\}$ and
 let $b= \inf \{t: J(t) \in \partial G_{n}
\hspace{.1in}\mbox{and}\hspace{.1in} t\geq s\}$.
 Then $J_{n} = J([a,b])$
is a Jordan arc in $\overline G_{n}$ with endpoints in $\partial
G_{n}$. Since $\partial G_{n}$ is an analytic curve, it is evident
that  $J_{n}$ separates $G_{n}$ into two disjoint domains. Let $x$
and $y$ be points from them, respectively, and let $\gamma$ be a
path which joins $x$ and $y$ and is contained in  $  G_{n}$ Since
$x$ and $y$ are separated by two domains enclosed by $\partial G_{n}
\cap J_{n}$, it follows that $\gamma \cap J_{n} \neq \emptyset$.
Hence, $\gamma$ is not contained in $G-J$ and consequently $G-J$ is
not connected. Lastly, let $G-J =G_{1}\cup G_{2}$, where $G_{i}$ are
disjoint domains. Since   {\bf C} - $G_{1}$ =  ({\bf C} $- G $)
$\cup J  \cup G_{2}$ is connected, it follows that  $G_{1}$ is
simply connected.
Similarly, $G_{2}$ is also simply connected.\\

The proof above is taken from \cite{carl-c}.

 \noindent   {\bf  B}). \  {\em
Let $f$ be an analytic function on $G$ and let $E\subset\partial G$
 such that $\partial G-E$ is a Jordan arc. If   there is a constant
$c$ such that for each $\lambda\in E$, $f(z)\rightarrow c$ as $z$
approaches to $\lambda$ from the inside of $G$, then $f$ is a
constant
function on $G$.}\\

{\bf Proof}. Firstly, if $G$ is a disk, then the conclusion follows
easily from the reflection principle.  For a general $G$, observe
that we can find an arc $J$  which is  a portion of a circle
   such that $J$ is a crosscut  of $G$ and its endpoints are in $E$.  $G-J$ has a component $U$ for which
    $\partial U\cap E$   connected. Let $g$ be a conformal map from $D$ onto $U$.
Then  $g^{-1}$ can extend continuously to $J^{\circ} $.\ Let
$L=g^{-1}(J^{\circ})$, then it is a subarc of $\partial D$. Now the
hypothesis implies that  for each $b\in
\partial D-\overline L$, $f\circ g(z)\rightarrow c$ as
$z\rightarrow b$ from the inside of $\varphi(U)$, and this  infers
that $f\circ  g=c$. Hence $f=c$ on $G$.\\

 \noindent {\bf  C}).\  {\em Let $I$ be a subarc of $\partial D$, then there is $\lambda\in I^{\circ}$
such that $\psi$ has radical limit at $\lambda$.}\\

{\bf Proof}.  Suppose the contrary. For $b\in I$ and let
 $J_{b}$ be the radius ending at $b$, then
  $ \overline {\psi(J_{b})}\cap \partial G $ is connected.
  Let $\Gamma_{b}=  \overline {\psi(J_{b})}\cap \partial G$.
  Note,     for  $n\geq 1$, the number of balls of radius great than $\frac{1}{n}$ and mutually disjoint  in $
  \overline G$ is finite.
 So,  there are   $b_{1}$ and $b_{2}$ in $I$, such that $\Gamma_{b_{1}}\cap \Gamma_{b_{2}}\neq
 \emptyset$. Notice that $J_{b_{1}}\cup J_{b_{2}}$ is a crosscut
 that separates $D$ into two parts, and let $W$ be the part whose
 boundary contains an entire arc which joins $b_{1} $ and $b_{2}$.
 Let $a\in \Gamma_{b_{1}}\cap\Gamma_{b_{2}}$ and set
 $g(z)=\frac{1}{\psi(z)-a}$. $g$ maps $W$ onto a simply connected
 domain $\Omega$, and $\partial \Omega= g(J_{b_{1}})\cap
 g(J_{b_{2}})\cup \{\infty\}$. Let $h$ be the Riemann map that
 maps $\Omega$ onto $D$, and set  $\alpha=h\circ g$, then   $\partial D= h(J_{b_{1}})\cap
 h(J_{b_{2}})\cup \{h(\infty)\}$.
 So we have that  $ lim_{z\rightarrow w}\alpha(z)=h(\infty)$ for
 $w\in \partial W\cap\partial D$ and it follows by Appendix B) that $\alpha$ must be a
 constant.
This is a contradiction.\\

 \noindent
{\bf  D}). \  {\em If $\lambda\in\partial_{a}G$ and   $J$ is a is a
Jordan arc such that $J\subset\overline G$ and $J\cap \partial
G=\{\lambda\}$, then $\lim_{z\rightarrow \lambda}\varphi(z)$ exists,
where the limit is taken\ along with $J$.}\\

{\bf Proof}.\ Suppose  the limit does not exist. Let $I=\overline
{\varphi(J^{\circ})}\cap\partial D$, then $I$ is a subarc. By virtue
of Appendix C, we can find a crosscut $L$ such that its endpoints
are on $I^{\circ}$ and $\psi  $ maps $L$ onto an open arc in $G$.
Let $\gamma=\overline{\psi(L )}$.  Then it is not difficult to see
that $\gamma$ is a Jordan curve and $\gamma\cap\partial
G=\{\lambda\}$.
Let $U$  be the domain enclosed by $\gamma$.
 Then $\partial U=\gamma$, and this implies that $\psi=\lambda$ on $\varphi(\gamma)-L$.
So it follows that $\varphi$ is a constant function, a
contradiction.\\

  \noindent {\bf  E}).  {\em  If $J$ is a crosscut of $G$, then $ \varphi(J) $ is a
crosscut of $D$.}\\

 {\bf Proof}.\
 Let $a$ and $b$ be the end points of $J$. If the lemma is not
true, then $\overline{\varphi(J)}$ is a Jordan curve
  and $\overline{\varphi(J)}\cap\partial D =\{\lambda\}$ for some
  $\lambda\in\partial D$.
Let $W$ be the domain enclosed by $\overline{\varphi(J)}$ and let
$U=\varphi^{-1}(W)$. Then $\varphi$ is constantly equal to $\lambda$
on $\partial U\cap \partial G$, now it follows from Appendix B) that
$\varphi$ is a constant. This is a contradiction.\\

{\noindent \bf Remark.} The author would like to point out that  all
of our proofs in this article do not involve the Lebesgue
integration theory. According to both \cite{gray} and \cite{sh},
Carathe\"{o}dory thought it was not possible for Osgood to prove his
conjecture at the time he gave since the Lebesgue   theory did not
exist yet.

\end{document}